\documentclass[a4paper,10pt]{article}

\makeatletter

\usepackage[utf8]{inputenc}
\usepackage[T1]{fontenc}
\usepackage[czech,english]{babel}

\usepackage{geometry}
\usepackage{ifpdf}
\ifpdf
        \@ifpackageloaded{hyperref}{}{\usepackage[pdftex]{hyperref}}
	\@ifpackageloaded{xcolor}{}{\usepackage[pdftex,usenames,svgnames,table]{xcolor}}
        \@ifpackageloaded{graphicx}{}{\usepackage[pdftex]{graphicx}}
	\usepackage{epstopdf}
\else
        \@ifpackageloaded{hyperref}{}{\usepackage[hypertex]{hyperref}}
	\@ifpackageloaded{xcolor}{}{\usepackage[dvips,usenames,svgnames,table]{xcolor}}
        \@ifpackageloaded{graphicx}{}{\usepackage[dvips]{graphicx}}
\fi
\graphicspath{{fig/}}

\definecolor{blue}{rgb}{0,0,0.5} 
\hypersetup{colorlinks=true,allcolors=blue}
\usepackage{xurl} 
\usepackage{subcaption}

\usepackage{amssymb,amsfonts}
\usepackage{lmodern}
\usepackage{dsfont}

\usepackage{fancyhdr}

\usepackage{amsmath,amsthm,amscd,upref}
\usepackage{mathtools}
\mathtoolsset{showonlyrefs=true} 

\usepackage{stackengine} 
\stackMath
\newcommand\ntilde[2][2]{%
 \def\useanchorwidth{T}%
  \ifnum#1>1%
    \stackon[-1.3ex]{\ntilde[\numexpr#1-1\relax]{#2}}{\mathchar"307E\kern-.5pt}%
  \else%
    \stackon[-1ex]{#2}{\mathchar"307E\kern-.5pt}%
  \fi%
}
\let\ubar=\underline

\usepackage{paralist} 
\usepackage{enumitem}

\usepackage{tikz}

\usepackage[round,sort]{natbib}

\newcommand{\doi}[1]{DOI~\href{\detokenize{https://doi.org/#1}}{\detokenize{#1}}}
\newcommand{\zblnumber}[1]{Zbl~\href{\detokenize{https://zbmath.org/?q=an:#1}}{\detokenize{#1}}}
\newcommand{\mrnumber}[1]{\href{\detokenize{https://www.ams.org/mathscinet-getitem?mr=#1}}{\detokenize{MR#1}}}

\definecolor{blue}{rgb}{0.1,0.1,0.5}
\definecolor{favorange}{cmyk}{0.02,0.22,1.0,0.15}
\definecolor{lightyellow}{cmyk}{0,0,0.1,0}

\definecolor{myfullred}{RGB}{244,63,43}
\definecolor{myhalfred}{RGB}{158,22,7}
\definecolor{myfullorange}{RGB}{255,128,0}
\definecolor{myhalforange}{RGB}{147,73,0}
\definecolor{mygreen}{RGB}{70,180,5}
\definecolor{mylilas}{RGB}{170,55,241}

\colorlet{pagebgcolor}{white}
\colorlet{textcolor}{black}
\colorlet{titlebgcolor}{pagebgcolor}

\usepackage{xspace}

\usepackage{soulutf8}
\sethlcolor{LimeGreen}

\soulregister\cite7%
\soulregister\citep7%
\soulregister\eqref7%
\soulregister\pageref7%
\soulregister\ref7%

\@ifundefined{ccode}{%
\newcommand{\ccode}[2]{\par
        \vspace*{8pt}
        {{\leftskip18pt\rightskip\leftskip
        \noindent{\it #1}\/: #2\par}}\par}
}{}
\@ifundefined{keywords}{\newcommand{\keywords}[1]{\ccode{Keywords}{#1}}}{}
\@ifundefined{email}{\newcommand{\email}[1]{\href{mailto:#1}{#1}}}{}

\@ifclassloaded{article}{%

\geometry{a4paper,top=20mm,left=20mm,right=20mm,bottom=20mm} 
\setlength{\headheight}{18.0pt} 

\rhead[]{\fancyplain{}{\slshape\rightmark}}
\chead{}
\lhead[\fancyplain{}{\slshape\leftmark}]{}
\cfoot{\fancyplain{\bfseries --~\thepage~--}{\bfseries --~\thepage~--}}

\usepackage{authblk}

\renewcommand\section{\@startsection {section}{1}{\z@}%
                                   {-3.5ex \@plus -1ex \@minus -.2ex}%
                                   {2.3ex \@plus.2ex}%
                                   {\color{blue}\normalfont\Large\bfseries }}

}{}

\@ifclassloaded{beamer}{%
\usetheme{default}
\setbeamertemplate{navigation symbols}{}
\setbeamertemplate{sidebar left}{\rlap{\hskip1mm\insertlogo}}
\setbeamertemplate{sidebar right}{}
\newcommand{\myfootlinefull}
{%
  \leavevmode%
  \hbox{\begin{beamercolorbox}[wd=\paperwidth,ht=2.5ex,dp=1.125ex,leftskip=3mm,rightskip=3mm]{author in head/foot}%
  \usebeamerfont{author in head/foot}\insertshortdate\hfill\insertshortauthor:\:
  \usebeamerfont{title in head/foot}\insertshorttitle\hfill \insertframenumber~/~\inserttotalframenumber
  \end{beamercolorbox}}%
  \vskip0pt%
}

\setbeamertemplate{footline}{\myfootlinefull}
\setbeamercolor{background canvas}{bg=pagebgcolor}
\setbeamercolor{alerted text}{fg=Navy}
\setbeamercolor{structure}{fg=Navy}
\setbeamercolor{frametitle}{fg=textcolor}
\setbeamercolor{author in head/foot}{fg=textcolor}
\setbeamertemplate{frametitle}[default][right,rightskip=1mm]
\addtobeamertemplate{frametitle}{}{\noindent\makebox[\linewidth]{\textcolor{favorange}{\rule{\paperwidth}{1pt}}}}
\setbeamerfont{title}{shape=\scshape}

\setbeamertemplate{bibliography item}{}
\setbeamersize{text margin left=5mm,text margin right=5mm}

\addtobeamertemplate{theorem begin}{%
    \setbeamercolor{block title}{fg=black,bg=favorange}%
    \setbeamercolor{block body}{fg=black,bg=lightyellow}%
}{}

\logo{%
\begin{tabular}{p{26mm}p{30mm}}%
\vspace{0mm}\href{http://www.fav.zcu.cz/en}{\includegraphics[height=13mm,trim=5mm 5mm 7mm 5mm,clip]{logo-fav_en}}&%
\vspace{1.5mm}\href{http://www.kma.zcu.cz/en}{\includegraphics[height=13mm,trim=0mm 0mm 135mm 0mm,clip]{logo-kma_en}}
\end{tabular}
}
\titlegraphic{%
\begin{tabular}{p{26mm}p{34mm}p{40mm}p{50mm}}%
\vspace{0mm}\href{http://www.fav.zcu.cz/en}{\includegraphics[height=13mm,trim=5mm 5mm 7mm 5mm,clip]{logo-fav_en}}&%
\vspace{1.5mm}\href{http://www.kma.zcu.cz/en}{\includegraphics[height=13mm,trim=0mm 0mm 135mm 0mm,clip]{logo-kma_en}}&%
\vspace{0mm}\href{http://www.ntis.zcu.cz/en}{\includegraphics[height=12mm,trim=0mm 0mm 0mm 0mm,clip]{logo-ntis}}&%
\vspace{1.5mm}\href{http://www.uni-rostock.de/en}{\includegraphics[height=10mm,trim=0mm 0mm 80mm 0mm,clip]{logo-rostock}}
\end{tabular}
}

}{} 

\IfFileExists{definitions.tex}{\input{definitions}}{}

\makeatother

\newcommand{\jpTitle}{Monotone iteration scheme for nonlinear PDEs in risk models}

\newcommand{\jpKeywords}{value adjustment; nonlinear PDE; monotone iterations method; subsolution; supersolution}
\newcommand{\jpMSC}{35A16; 35K58; 91G40; 91G20}
\newcommand{\jpJEL}{G12; C63}

\newcommand{\jpDate}{29 June 2023}

\author[1]{Falko Baustian} 
\author[2]{Jan Posp\'{\i}\v{s}il\thanks{Corresponding author, \email{honik@kma.zcu.cz}}} 
\author[2]{Vladim\'{\i}r \v{S}v\'{\i}gler}
\affil[1]{Institute of Mathematics, University of Rostock, Ulmenstra\ss e 69, 18057 Rostock, Germany}
\affil[2]{NTIS - New Technologies for the Information Society, Faculty of Applied Sciences, \authorcr University of West Bohemia, Univerzitn\'{\i} 2732/8, 301 00 Plze\v{n}, Czech Republic,\vspace*{3pt}}

\title{\textcolor{Navy}{\textsc{\jpTitle}}}
\date{\jpDate}

\begin{document}

\maketitle

\begin{center}
\end{center}

\begin{abstract}
In this paper we study nonlinear partial differential equations (PDEs) that are used to model different value adjustments denoted generally as xVA. These adjustments are nowadays commonly added to the risk-free financial derivative values and the PDE approach allows their easy incorporation. The aim of this paper is to apply the method of monotone iterations with sub- and supersolutions to the nonlinear Black-Scholes-type equation that occurs especially in the counterparty risk models. We introduce a monotone iteration scheme with semi-explicit solution formulas for each iteration step. Moreover, we show that the problem greatly simplifies for contracts with non-negative payoffs. To show the viability of the approach we apply our method to the call option, the forward, and the gap option.

\end{abstract}

\keywords{\jpKeywords}
\ccode{MSC classification}{\jpMSC}
\ccode{JEL classification}{\jpJEL}

\setcounter{tocdepth}{2}
\tableofcontents
\clearpage


\section{Introduction}\label{sec:introduction}

The \emph{method of sub- and supersolutions} is a technique used to obtain approximate solutions for differential equations. It is particularly useful when it is difficult or impossible to find an exact solution to a given equation. The basic idea is to construct two functions, a subsolution and a supersolution, that sandwich the true solution of the equation. By comparing these functions, we can gain insights into the behaviour of the solution and estimates for it. The method of sub- and supersolutions was known for a long time in the context of ordinary differential equations, see for example the historical overview in Chapter 1 of the book by \cite{DeCosterHabets2006}. The first author who used it in partial differential equations (PDEs) was to the best to our knowledge \cite{Nagumo1954}. Some years later, especially after the seminal paper by \cite{KellerCohen1967} was published, the method became a standard tool. One of the first papers that focused on monotonicity of super- and subsolutions in nonlinear elliptic and parabolic boundary value problems (BVPs) was written by \cite{Sattinger1972}. In this paper, he among others presented as an open problem the question of existence of a solution in presence of sub- and supersolutions that do not satisfy the natural well-ordering relation.

Independently of the technique of sub- and supersolutions, the \emph{fixed point method}, became popular at the late 1970s, see e.g. the book by \cite{KinderlehrerStampacchia1980} or \citep[reprint]{KinderlehrerStampacchia2000}. In the fixed point method, the idea is to iteratively construct a sequence of approximations that converge to the solution of the PDE. A general outline of the fixed point method can be written as:

\begin{enumerate}
\item Start with an initial guess or approximation for the solution, denoted by $u_0$.   
\item Perform a sequence of iterations to improve the approximation:
\begin{itemize}
\item At each iteration $k$, solve a modified PDE, denoted by $\mathcal{L}(u_{k+1}) = f(u_k)$, where $\mathcal{L}$ is the differential operator and $f$ is the corresponding reaction term.  
\item The modified PDE incorporates the current approximation $u_k$ and aims to correct its deviation from the true solution. 
\item The solution to the modified PDE provides the updated approximation $u_{k+1}$ for the next iteration.
\end{itemize}
\item Continue the iterations until the sequence converges, i.e., until the difference between successive approximations becomes sufficiently small.
\end{enumerate}

The key aspect of the fixed point method is that each iteration produces a new approximation that improves upon the previous one. The sequence of approximations generally converges to the true solution if certain conditions, such as monotonicity, coercivity, or compactness, are satisfied. The fixed point method is particularly useful for PDEs with monotone operators, as it exploits the order-preserving nature of these operators to ensure convergence. However, it may not be applicable or efficient for all types of PDEs, especially those with non-monotone or highly nonlinear operators. Despite some inconsistency of the method names in the literature, we call the fixed point method for monotone operators the \emph{monotone method}. 

From the numerical point of view, monotone methods rely on the weak comparison principle  
that ensures convergence and stability of the solution. 
Monotonicity preserving in numerical methods for PDEs is discussed for finite differences methods in the book by \cite{Strikwerda2004} and for finite element methods in the book by \cite{Thomee2006}. We also refer to the monograph by \cite{Pao92}. Numerical methods for nonlinear PDEs in finance were also studied by \cite{Forsyth2012}.

While Kinderlehrer and Stampacchia's contribution was significant in the development of the theory and applications of monotone iterations, it is important to note that the method itself has evolved through the contributions of many researchers over time. Monotone methods were applied to various PDEs that appear in multiple disciplines, for example to the 
porous medium equation in hydrology \citep{PeletierTroy2001},
Fisher-Kolmogorov–Petrovsky–Piskunov equation in population genetics \citep{Gomez2011},
Einstein-Lichnerowicz equation in general relativity \citep{Premoselli2014}, 
Allen-Cahn equation in material engineering \citep{ChenGuoNinomiyaYao2018}, 
etc.

It is worth to mention, that the monotone method does not have to deal with sub- and supersolutions in general at all. Instead, it focuses on iteratively improving an initial approximation to converge to the solution of the PDE. However, our aim is to consider the \emph{method of monotone iterations} (or monotone sequences) \emph{with sub- and supersolutions} whose idea is described for example in \cite[Chap. V]{DeCosterHabets2006} that summarizes the earlier results \cite{CherpionDeCosterHabets1999,CherpionDeCosterHabets2001}. In this method, an increasing sequence of subsolutions and a decreasing sequence of supersolutions are constructed similarly as in the fixed point method algorithm above. Due to the boundedness of both sequences, they both have a limit and converge to a fixed point of $\mathcal{L}-f$. 

Non-linear PDEs in finance appear in a number of related but sometimes also distinct applications. These applications include
pricing derivatives contracts with a one-way Credit Support Annex (CSA) \citep{Piterbarg10},
risky closeouts for Credit Value Adjustments (CVAs) \citep{BurgardKjaer11PDE}, 
derivatives pricing with differential borrowing and lending rates \citep{Mercurio13}, 
and accounting-consistent valuation \citep{BurgardKjaer15}.
In this paper, we study nonlinear PDE of the Black-Scholes type that can be used to model value adjustments. Different value adjustments denoted generally as xVA are nowadays added to the risk-free financial derivative values to take the complex market behaviour into consideration and the PDE approach allows their easy incorporation. 

The aim of this paper is to apply the method of monotone iterations with sub- and supersolutions to the nonlinear Black-Scholes-type equation that occurs in the counterparty risk models. Motivated by the work of \cite{AlziaryTakac22} we provide a new insights to the PDEs studied by \cite[Sec. 9]{BaustianTakac21} and \cite{BaustianFenclPospisilSvigler22wilmott}. We also refer to the works by \cite{BurgardKjaer11}, \cite[Section~3]{BurgardKjaer11PDE}, and \cite{BurgardKjaer11PDE-add} for details concerning xVA modelling and to \cite{Arregui17} and \cite{BaustianFenclPospisilSvigler22wilmott} for numerical results. 

The paper is structured in the following way. We introduce the PDE model in Section~\ref{sec:preliminaries} together with the solution for the linear case studied in \citep{BaustianFenclPospisilSvigler22wilmott}. In Section~\ref{sec:methodology}, we introduce the new monotone iteration scheme for the nonlinear problem and provide formulas for selected derivatives contracts. Numerical results are presented in Section~\ref{sec:results}. A conclusion is given in Section~\ref{sec:conclusion}.

\section{Preliminaries}\label{sec:preliminaries}

We consider a derivative contract between an issuer $B$ and a counterparty $C$ in the setting of the Black-Scholes model but incorporate funding costs and the bilateral counterparty risk. The market consists of a risk-free bond with a constant interest rate $r>0$ and a risky asset $S$, e.g., a stock, that satifies the stochastic differential equation
$$
\mathrm{d}S_t = (q_S-\gamma_S)S_t\,\mathrm{d}t + \sigma S_t\,\mathrm{d}W_t,
$$
where the effective financing rate $q_S-\gamma_S$ and the volatility $\sigma>0$ are constant and $W_t$ is a standard Wiener process. A European-style derivative is a contract whose payoff depends on the value of the underlying asset $S$ at some fixed maturity time $T$. The arbitrage-free or fair price $V(t,S)$ of such a derivative in the classical Black-Scholes model without funding costs or the possibility of default satisfies the Black-Scholes equation
\begin{equation*}
    \frac{\partial}{\partial t}V+\mathcal{A}V-rV=0
\end{equation*}
in $(t,S)\in(0,T)\times(0,\infty)$ with the final data $V(T,S)=H(S)$ given by the payoff function $H$. Here, the differential operator $\mathcal{A}$ is the Black-Scholes operator
\begin{equation*}
    \mathcal{A}=\frac12\sigma^2S^2\frac{\partial^2}{\partial S^2}+(q_S-\gamma_S)S\frac{\partial}{\partial S}.
\end{equation*}
We study the adjusted price $U(t,S)$ of a derivative contract under funding costs and the bilateral counterparty risk. According to \cite{BurgardKjaer11PDE}, the adjusted price $U$ satisfies the PDE
\begin{equation*}
    \frac{\partial}{\partial t}U+\mathcal{A}U-(r+\lambda_B+\lambda_C)U=(R_B\lambda_B+\lambda_C)M^--(\lambda_B+R_C\lambda_C)M^++s_FM^+ 
\end{equation*}
in $(t,S)\in(0,T)\times(0,\infty)$ with the same final value $U(T,S)=H(S)$, where $\lambda_B,\lambda_C$ are the default intensities of the counterparties, $R_B,R_C\in[0,1]$ are the recovery rates of the derivative positions, and $s_F$ is the funding spread. The term $M$ describes the mark-to-market value $M$ of the derivative at default. In this article, we use the standard abbreviations $\varphi^+=\max(\varphi,0)$ and $\varphi^-=\max(-\varphi,0)$ for positive and negative parts of a function $\varphi$ that slightly differs from the notation in \cite{BurgardKjaer11PDE}. 
A common choice for the mark-to-market value is $M=V$, e.g., for contracts under the 2002 International Swaps and Derivatives Association (ISDA) Master Agreement, and it yields the linear equation
\begin{equation}\label{eq:V_lin}
    \frac{\partial}{\partial t}U+\mathcal{A}U-(r+\lambda_B+\lambda_C)U=(R_B\lambda_B+\lambda_C)V^--(\lambda_B+R_C\lambda_C)V^++s_FV^+. 
\end{equation}
In our previous work \cite{BaustianFenclPospisilSvigler22wilmott}, we derived a solution formula for the difference $U-V$ in this linear case and compared it to standard pricing methods for adjusted contracts. Another reasonable choice for the mark-to-market value, $M=U$, gives the semilinear PDE
\begin{equation}\label{eq:V_nonlin}
    \frac{\partial}{\partial t}U+\mathcal{A}U-rU=-(1-R_B)\lambda_BU^-+(1-R_C)\lambda_CU^++s_FU^+. 
\end{equation}
This equation was recently studied in \cite{AlziaryTakac22} where the authors transform \eqref{eq:V_nonlin} into a parabolic problem with a monotone nonlinear reaction function. In \cite[Sec. 3.2]{AlziaryTakac22}, they construct an iteration scheme of monotone sequences of \emph{subsolutions} and \emph{supersolutions} which convergence monotonically to the solution of the nonlinear problem. We call $\ubar{u}$ the subsolution of an initial value problem for a differential equation $\mathcal{L}u=f(u)$ with a differential operator $\mathcal{L}$, a reaction function $f$, and the initial data $u(0)=u_0$, when it satisfies $\mathcal{L}\ubar{u}\leq f(\ubar{u})$ and $\ubar{u}(0)\leq u_0$. A supersolution $\bar{u}$ of the same problem satisfies $\mathcal{L}\bar{u}\geq f(\bar{u})$ and $\bar{u}(0)\geq u_0$.

We give a short overview of the contracts that we will study. Many articles on derivative pricing under value adjustments use European call options to validate their methods. European call options give the issuer the right but not the obligation to buy an asset at the maturity time $T$ of the contract for a fixed strike price $K$. The payoff function of the call option is
$$
H_{\mathrm{call}}(S)=(S-K)^+
$$
and its risk-free price is given by the famous Black-Scholes formula
$$
V_{\mathrm{call}}(t,S)=\mathrm{e}^{(q_s-\gamma_S-r)(T-t)}S\mathcal{N}(d_1)-\mathrm{e}^{-r(T-t)}K\mathcal{N}(d_2),
$$
where $\mathcal{N}$ is the standard normal cumulative distribution function, and
$$
d_{1,2}=\frac{\ln(S/K)+(r\pm\sigma^2/2)(T-t)}{\sigma\sqrt{T-t}}.
$$
The European call option and most other derivatives have a non-negative payoff. In Section \ref{ssec:nonneg}, we show that the adjusted price for such contracts can be derived easily from the risk-free price. Hence, we are mostly interested in contracts with sign-changing payoffs. Forwards give the issuer the obligation to buy an asset at the fixed future time $T$ for a pre-agreed forward price $K$. The payoff function of a forward is
$$
H_{\mathrm{fwd}}(S)=S-K
$$
and its price in the risk-free setting is given by
$$
V_{\mathrm{fwd}}(t,S)=\mathrm{e}^{(q_s-\gamma_S-r)(T-t)}S-\mathrm{e}^{-r(T-t)}K.
$$
A gap call option with a strike price $K_s$ which is higher than the trigger price $K_t$ has also a sign-changing payoff
$$
H_{\mathrm{gap}}(S)=(S-K_s)\cdot\mathds{1}_{\{S\geq K_t\}}.
$$
The risk-free price of the gap is given by
$$
V_{\mathrm{gap}}(t,S)=\mathrm{e}^{(q_s-\gamma_S-r)(T-t)}S\mathcal{N}(d_1)-\mathrm{e}^{-r(T-t)}K_s\mathcal{N}(d_2),
$$
where $d_1$ and $d_2$ are defined as in the Black-Scholes formula but depend on $K_t$ instead of $K$. We remark that the payoffs of the call option, the forward, and the gap call option satisfy the at most linear growth in infinity and local boundedness at zero, see also hypothesis $(v_0)$ from \cite{AlziaryTakac22}.

\subsection{Solution formula for the linear problem}\label{ssec:lin}

In \cite{BaustianFenclPospisilSvigler22wilmott}, we established a solution formula for a linear problem similar to the PDE \eqref{eq:V_lin}. The formula is based on a transformation to the heat equation. We can re-utilize this approach, see \cite[Sec. 3.1]{BaustianFenclPospisilSvigler22wilmott}, to derive a solution formula for the linear case.

We begin by introducing the new variables $x=\ln S$ and $\tau=T-t$ and the new function $u(\tau,x)=U(t,S)$ to transform the linear PDE \eqref{eq:V_lin} to the equation
\begin{equation*}
    \frac{\partial}{\partial\tau}u - \frac12\sigma^2\frac{\partial^2}{\partial x^2}u - \rho\frac{\partial}{\partial x}u + (r+\lambda_B+\lambda_C)u = f(V(T-\tau,\mathrm{e}^x))
\end{equation*}
with the abbreviation $\rho=q_S-\gamma_S-\frac12\sigma^2$ and the inhomogeneity
$$
f(V) = -(R_B\lambda_B+\lambda_C)V^- + (\lambda_B+R_C\lambda_C)V^+ - s_FV^+.
$$
The function $u$ satisfies the initial condition $u(0,x)=h(x)$ with the transformed payoff $h(x)=H(S)$.
Next, we set $\tilde{u}(\tau,x)=\mathrm{e}^{(r+\lambda_B+\lambda_C)\tau}u(\tau,x)$ and obtain
\begin{equation*}
    \frac{\partial}{\partial\tau}\tilde{u} - \frac12\sigma^2\frac{\partial^2}{\partial x^2}\tilde{u} - \rho\frac{\partial}{\partial x}\tilde{u} = f(V(T-\tau,\mathrm{e}^x))\mathrm{e}^{(r+\lambda_B+\lambda_C)\tau}.
\end{equation*}
Finally, we use the shift $\tilde{x} = x + \rho \tau$, i.e., $\check{u}(\tau,\tilde{x}) = \check{u}(\tau,x+\rho \tau)=  \tilde{u}(\tau,x)$. This leads to 
$$
\frac{\partial}{\partial \tau} \check{u}(\tau,\tilde{x}) = \frac{\partial}{\partial \tau} \tilde{u}(\tau,x) - \rho \frac{\partial}{\partial x} \tilde{u}(\tau,x). 
$$
Formally omitting the tilde sign at $x$ results in the initial value problem for the inhomogeneous heat equation
\begin{equation*}
    \frac{\partial}{\partial\tau}\check{u} - \frac12\sigma^2\frac{\partial^2}{\partial x^2}\check{u} = f(V(T-\tau,\mathrm{e}^x))\mathrm{e}^{(r+\lambda_B+\lambda_C)\tau}
\end{equation*}
with the initial data $\check{u}(0,x)=\tilde{u}(0,x)=u(0,x)=h(x)$. The solution to this problem is given by Duhamel's formula 
\begin{align}\label{eq:u_lin}
\begin{split}
    \check{u}(\tau,x) &= \int_{-\infty}^{\infty}\frac{1}{\sqrt{2\pi\sigma^2\tau}}\exp\left(-\frac{(x-y)^2}{2\sigma^2\tau}\right)h(y)\,\mathrm{d}y \\
    &+ \int_0^\tau\int_{-\infty}^{\infty}\frac{1}{\sqrt{2\pi\sigma^2(\tau-s)}}\exp\left(-\frac{(x-y)^2}{2\sigma^2(\tau-s)}\right)f(V(T-s,\mathrm{e}^y))\mathrm{e}^{(r+\lambda_B+\lambda_C)s}\,\mathrm{d}y\,\mathrm{d}s.
\end{split}
\end{align}
We obtain the following formula for $u$
\begin{align}\label{e:lin}
\begin{split}
    u(\tau,x) &= \frac{1}{\exp(-(r+\lambda_B+\lambda_C)\tau)}\int_{-\infty}^{\infty}\frac{1}{\sqrt{2\pi\sigma^2\tau}}\exp\left(-\frac{(x+\rho\tau-y)^2}{2\sigma^2\tau}\right)h(y)\,\mathrm{d}y \\
    &+ \int_0^\tau\int_{-\infty}^{\infty}\frac{1}{\sqrt{2\pi\sigma^2(\tau-s)}}\exp\left(-\frac{(x+\rho\tau-y)^2}{2\sigma^2(\tau-s)}\right)f(V(T-s,\mathrm{e}^y))\mathrm{e}^{-(r+\lambda_B+\lambda_C)(\tau-s)}\,\mathrm{d}y\,\mathrm{d}s.
\end{split}
\end{align}

\section{Methodology}\label{sec:methodology}

We provide all information needed to implement the monotone method for nonlinear PDEs in risk models in this section. In Section \ref{ssec:nonlin}, we introduce our monotone iteration scheme. We discuss the problem for derivatives with non-negative payoffs and give easy solution formulas for the adjusted price of such contracts in Section \ref{ssec:nonneg}. In Section \ref{ssec:sc}, we consider specific contracts and derive estimates for the truncation error. Numerical aspects of the iteration scheme are discussed in Section \ref{ssec:num}.

\subsection{Monotone iteration scheme for the nonlinear problem} \label{ssec:nonlin}

We solve the nonlinear problem \eqref{eq:V_nonlin} by establishing a monotone convergence scheme based on the results in \cite{AlziaryTakac22}. We have to transform the PDE \eqref{eq:V_nonlin} to an equation with a monotone reaction function. For that purpose, we introduce $u(\tau,x)=U(t,S)$ and add the term $c_Mu$ with a constant $c_M>0$, satisfying $c_M\geq\max\left\{ (1-R_B)\lambda_B, (1-R_C)\lambda_C + s_F \right\}$, to both sides of equation \eqref{eq:V_nonlin} to obtain the problem
\begin{equation}\label{e:check-v}
\begin{cases}
\dfrac{\partial}{\partial \tau}u + \mathcal{A}u = g(u), \\
u(0,x)=h(x)=H(\mathrm{e}^x),
\end{cases}
\end{equation}
in which 
\begin{equation*}
\mathcal{A}u =  - \frac12\sigma^2\frac{\partial^2}{\partial x^2}u - \rho\frac{\partial}{\partial x}u + (r+c_M)u \\
\end{equation*}
and
\begin{equation}\label{e:g}
    g(u) = (1-R_B)\lambda_Bu^- - (1-R_C)\lambda_Cu^+ - s_Fu^+ + c_Mu.
\end{equation}
The added term $c_Mu$ guarantees that the inhomogeneity $g$ is a monotone increasing mapping, see \cite[Section 2]{AlziaryTakac22}, where the constant $c_M$ is chosen exactly as $c_M=\max\left\{ (1-R_B)\lambda_B, (1-R_C)\lambda_C + s_F \right\}$.

We say that some function $\underline{u}$ is a subsolution to~\eqref{e:check-v} if 
\begin{enumerate}[label={(\roman*)}]
    \item it is dominated by the initial condition $\underline{u}(0,x) \le h(x)$ for every $x\in\mathbb{R}$, 
    \item the solution residual 
    \begin{equation}\label{e:residual}
    \mathcal{D}(\underline{u}):=\dfrac{\partial}{\partial \tau} \underline{u} + \mathcal{A} \underline{u} - g(\underline{u})
    \end{equation}
    is nonpositive, i.e., $\mathcal{D}(\underline{u})\le 0$. 
\end{enumerate}
Analogously, a function $\overline{u}$ is a supersolution to~\eqref{e:check-v} if $\overline{u}(0,x) \ge h(x)$ for every $x\in\mathbb{R}$ and $\mathcal{D}(\overline{u}) \ge 0$. Note that the considered financial contract provides the initial condition of~\eqref{e:check-v} and thus the initial sub- or supersolution must be constructed accordingly. 

Let $\underline{u_0}$ be the initial subsolution and let $\overline{u_0}$ be the initial supersolution. The results in~\cite{AlziaryTakac22} show that the successive solving of the problems 
\begin{equation}\label{e:sub}
\begin{cases}
\dfrac{\partial}{\partial \tau}\underline{u_{n+1}} + \mathcal{A}\underline{u_{n+1}} = g(\underline{u_{n}}), \quad n = 0, 1, \ldots \\
\underline{u_{n+1}}(0,x)=h(x)=H(\mathrm{e}^x),
\end{cases}
\end{equation}
results in the sequence of subsolutions $(\underline{u_n})_{n \in \mathbb{N}_0}$ satisfying $\mathcal{D}(\underline{u_n}) \le 0$ for all $n \in \mathbb{N}_0$ and 
$$
\underline{u_0} \le \underline{u_1} \le \underline{u_2} \ldots \le u^*,
$$ 
in which $u^*$ is the solution of~\eqref{e:check-v}, i.e., $\mathcal{D}(u^*)=0$. Analogously, solving  
\begin{equation}\label{e:super}
\begin{cases}
\dfrac{\partial}{\partial \tau}\overline{u_{n+1}} + \mathcal{A}\overline{u_{n+1}} = g(\overline{u_{n}}), \quad n = 0, 1, \ldots \\
\overline{u_{n+1}}(0,x)=h(x)=H(\mathrm{e}^x),
\end{cases}
\end{equation}
yields a sequence of supersolutions $(\overline{u_n})_{n \in \mathbb{N}_0}$ satisfying $\mathcal{D}(\overline{u_n}) \ge 0$ for all $n \in \mathbb{N}_0$ and 
$$
\overline{u_0} \ge \overline{u_1} \ge \overline{u_2} \ge \ldots \ge u^*.
$$ 
Moreover, $(\underline{u_n})_{n \in \mathbb{N}_0} \xrightarrow[n \to +\infty]{} u^*$ and $(\overline{u_n})_{n \in \mathbb{N}_0} \xrightarrow[n \to +\infty]{} u^*$. 

The problems~\eqref{e:sub} and~\eqref{e:super} are linear and we can thus use a approach similar to the one presented in Section~\ref{ssec:lin}. Without loss of generality, we define $u_n$ to be the $n$-th iteration of either sub- or supersolution. First, we introduce $\tilde{u}_{n+1}(\tau,x)=\mathrm{e}^{(r+c_M)\tau}u_{n+1}(\tau,x)$ and get
\begin{equation*}
\frac{\partial }{\partial \tau}\tilde{u}_{n+1} - \frac12\sigma^2\frac{\partial^2}{\partial x^2}\tilde{u}_{n+1} - \rho\frac{\partial}{\partial x}\tilde{u}_{n+1} 
= g(u_{n})e^{(r+c_M)\tau}.
\end{equation*}
Second, we reuse the shift $\check{u}_{n+1}(\tau,x+\rho \tau)=  \tilde{u}_{n+1}(\tau,x)$ and obtain the initial value problem
\begin{equation}
\label{eq:nonlin_tf}
\frac{\partial }{\partial \tau}\check{u}_{n+1} - \frac12\sigma^2\frac{\partial^2}{\partial x^2}\check{u}_{n+1}
= g(u_n(\tau,x-\rho\tau))e^{(r+c_M)\tau}
\end{equation}
with the initial data $\check{u}(0,x)= h(x)$. Next, we make use of Duhamel's formula
\begin{align*}
    \check{u}_{n+1}(\tau,x) &= \int_{-\infty}^{\infty}\frac{1}{\sqrt{2\pi\sigma^2\tau}}\exp\left(-\frac{(x-y)^2}{2\sigma^2\tau}\right)h(y)\,\mathrm{d}y \\
    &+ \int_0^\tau\int_{-\infty}^{\infty}\frac{1}{\sqrt{2\pi\sigma^2(\tau-s)}}\exp\left(-\frac{(x-y)^2}{2\sigma^2(\tau-s)}\right)g(u_n(s,y-\rho s)) \exp{((r + c_M)s)}\,\mathrm{d}y\,\mathrm{d}s, 
\end{align*}

Using the backward substitutions $u_{n+1}(\tau,x)= e^{-(r+c_M)\tau}\tilde{u}_{n+1}(\tau,x)=e^{-(r+c_M)\tau}\check{u}_{n+1}(\tau, x+ \rho \tau)$, we end up with the solution of~\eqref{e:sub}, resp.~\eqref{e:super}
\begin{align}
    u_{n+1}(\tau,x) &= \frac{1}{\exp((r+c_M)\tau)}\int_{-\infty}^{\infty}\frac{1}{\sqrt{2\pi\sigma^2\tau}}\exp\left(-\frac{(x+\rho\tau-y)^2}{2\sigma^2\tau}\right)h(y)\,\mathrm{d}y \\
    &+ \int_0^\tau\int_{-\infty}^{\infty}\frac{1}{\sqrt{2\pi\sigma^2(\tau-s)}}\exp\left(-\frac{(x+\rho\tau-y)^2}{2\sigma^2(\tau-s)}\right)g(u_n(s,y-\rho s)) \exp{(-(r + c_M)(\tau-s))}\,\mathrm{d}y\,\mathrm{d}s,  
\end{align}
in which we used the fact that the function $g(\, \cdot \,)$ is positive homogeneous of degree 1. We apply the substitution $z = y- \rho \tau$ to convert the iteration step into a form which is more suited for numerical computations
\begin{align}
\begin{split}\label{e:iter}
    u_{n+1}(\tau,x) &= \frac{1}{\exp((r+c_M)\tau)}\int_{-\infty}^{\infty}\frac{1}{\sqrt{2\pi\sigma^2\tau}}\exp\left(-\frac{(x+\rho\tau-y)^2}{2\sigma^2\tau}\right)h(y)\,\mathrm{d}y  \\
    &+ \int_0^\tau\int_{-\infty}^{\infty}\frac{1}{\sqrt{2\pi\sigma^2(\tau-s)}}\exp\left(-\frac{(x-z)^2}{2\sigma^2(\tau-s)}\right)g(u_n(s,z+\rho (\tau -s))) \exp{(-(r + c_M)(\tau-s))}\,\mathrm{d}z\,\mathrm{d}s.  
\end{split}
\end{align}

\subsection{Formulas for non-negative payoffs}\label{ssec:nonneg}

If we consider derivatives with non-negative payoffs then the PDEs \eqref{eq:V_lin} and \eqref{eq:V_nonlin} simplify significantly. First, we study the linear equation \eqref{eq:V_lin}. The risk-less price $V$ of a contract with a non-negative payoff is also non-negative due to the weak maximum principle. Hence, the inhomogeneity of \eqref{eq:V_lin} has the form $(s_F-\lambda_B-R_C\lambda_C)V$ and, using similar arguments as in \cite[Sec. 3.2]{BaustianFenclPospisilSvigler22wilmott}, the adjusted price $U$ is directly given by
\begin{equation*}
    U=\frac{((1-R_C)\lambda_C+s_F)\mathrm{e}^{-(\lambda_B+\lambda_C)(T-t)}+\lambda_B+R_C\lambda_C-s_F}{\lambda_B+\lambda_C}\, V.
\end{equation*}
Second, in the nonlinear case \eqref{eq:V_nonlin} the solution is due to \cite{AlziaryTakac22} the limit of a monotonically increasing sequence of subsolutions $\underline{U_0}\leq\underline{U_1}\leq\underline{U_2}\leq\dots$ and, thus, the solution must satisfy $U\geq\underline{U_0}$. Since we can choose $\underline{U_0}\equiv0$ as the inital subsolution of the converging sequence for any non-negative payoff function the solution is also non-negative. As a consequence, problem \eqref{eq:V_nonlin} becomes the linear homogeneous equation
\begin{equation*}
    \frac{\partial}{\partial t}U+\mathcal{A}U-(r+(1-R_C)\lambda_C+s_F)U=0 
\end{equation*}
with the adjusted price
\begin{equation*}
    U=\mathrm{e}^{-((1-R_C)\lambda_C+s_F)(T-t)} V.
\end{equation*}
as its solution.

\subsection{Specific contracts}\label{ssec:sc}

In this section, we give useful calculations for the three studied contracts: the call option, the forward, and the gap option. We also provide bounds for the truncation error for the second integral in \eqref{e:iter} in the monotone iteration method introduced in Section \ref{ssec:nonlin}.

For each contract we need an initial subsolution and an initial supersolution to run the monotone iteration scheme. A general example for such subsolution and supersolution is given in \cite[Example 3.2]{AlziaryTakac22}. Whereas, we are looking for particular functions suited for our regarded contracts. For the initial supersolution, we consider a function of the form $\overline{u_0}(\tau,x)=\exp(\overline{\lambda}\tau+x)$ with $\overline{\lambda}\in\mathbb{R}$. The respective solution residual
$$
\mathcal{D}(\overline{u_0})=(\overline{\lambda}-\frac12\sigma^2-\rho+r+(1-R_C)\lambda_C+s_F)\exp(\overline{\lambda}\tau+x)
$$
is non-negative for every $\overline{\lambda}\geq\frac12\sigma^2+\rho-r-(1-R_C)\lambda_C-s_F$ and we have $\overline{u_0}(0,x) = \mathrm{e}^x$ with
\begin{align*}
    \overline{u_0}(0,x) &\geq h_{\mathrm{call}}(x)=(\mathrm{e}^x-K)^+,\\
    \overline{u_0}(0,x) &\geq h_{\mathrm{fwd}}(x)=\mathrm{e}^x-K, \\
    \overline{u_0}(0,x) &\geq h_{\mathrm{gap}}(x)=(\mathrm{e}^x-K_s)\mathds{1}_{\{x\geq\ln K_t\}},
\end{align*}
i.e., $\overline{u_0}$ is a supersolution for all three contracts. An obvious choice for the initial subsolution for derivatives with non-negative payoff, as the call option, is $\underline{u_0}(\tau,x)\equiv0$ with $\mathcal{D}(\underline{u_0})=0$. For contracts with sign-changing payoffs, as the forward and the gap option, it is natural to assume that the payoff function is at least bounded from below by a constant. We consider a function of the form $\underline{u_0}(\tau,x)=-\underline{C}\exp(\underline{\lambda}\tau)$ with $\underline{C}\geq0$ as the initial subsolution. The subsolution residual
$$
\mathcal{D}(\underline{u_0})=-\underline{C}(\underline{\lambda}+r-(1-R_B)\lambda_B)\exp(\underline{\lambda}\tau)
$$
is non-positive for every $\underline{\lambda}\geq-r+(1-R_B)\lambda_B$ and also in the trivial case $\underline{C}=0$ corresponding to non-negative payoff functions. Hence, we obtain a subsolution for the forward with $\underline{C}=K$ and a subsolution for the gap option with $\underline{C}=K_s-K_t$, respectively.

We can make use of the initial subsolution and the initial supersolution to control the truncation error in the monotone iteration scheme, i.e., we can estimate the error that occurs when we calculate the second integral in \eqref{e:iter} only on a finite interval $[-L,L]$ with $L>0$. For this purpose, we suppose a priori bounds on the initial subsolution $\underline{u_0}$ and the initial supersolution $\overline{u_0}$ in accordance with the choices above. We assume
\begin{equation}\label{eq:bds}
-\underline{C}\mathrm{e}^{\underline{\lambda}\tau}\leq\underline{u_0}\leq\overline{u_0}\leq\overline{C}\mathrm{e}^{\overline{\lambda}\tau}\mathrm{e}^{ax}
\end{equation}
for some constants $\underline{C} \ge 0$, $\overline{C}>0$, and $\underline{\lambda},\overline{\lambda},a\in\mathbb{R}$. Under assumption \eqref{eq:bds}, the monotonicity of the iteration scheme guarantees
\begin{equation*}
-\underline{C}\mathrm{e}^{\underline{\lambda}\tau}\leq\underline{u_0}\leq\underline{u_n}\leq\overline{u_n}\leq\overline{u_0}\leq\overline{C}\mathrm{e}^{\overline{\lambda}\tau}\mathrm{e}^{ax} \quad\mbox{for every }n\in\mathbb{N}.
\end{equation*}
Let $u$ be a function satisfying $-\underline{C}\mathrm{e}^{\underline{\lambda}\tau}\leq u\leq\overline{C}\mathrm{e}^{\overline{\lambda}\tau}\mathrm{e}^{ax}$. The monotonicity of the inhomogeneity $g$ defined in \eqref{e:g} gives
$$
-\underline{C}(c_M-(1-R_B)\lambda_B)\mathrm{e}^{\underline{\lambda}\tau}\leq g(u)\leq \overline{C}(c_M-(1-R_c)\lambda_C-s_F)\mathrm{e}^{\overline{\lambda}\tau}\mathrm{e}^{ax}
$$
and we obtain the bound
$$
|g(u)|\leq C^+\mathrm{e}^{\overline{\lambda}\tau}\mathrm{e}^{ax} + C^-\mathrm{e}^{\underline{\lambda}\tau}
$$
with $C^+=\overline{C}(c_M-(1-R_c)\lambda_C-s_F)$ and $C^-=\underline{C}(c_M-(1-R_B)\lambda_B)$. We remark that for the particular choice $c_M=\max\left\{ (1-R_B)\lambda_B, (1-R_C)\lambda_C + s_F \right\}$ in \cite{AlziaryTakac22} rather $C^-$ or $C^+$ vanishes. We can now estimate the truncation error for the second integral in \eqref{e:iter} in each iteration step provided \eqref{eq:bds} holds. For better clearness of the display, we estimate the truncation error regarding $(L,\infty)$ and $(-\infty,L)$ separately. We have
\begin{align*}
    &\left|\int_0^\tau\int_{L}^{\infty}\frac{1}{\sqrt{2\pi\sigma^2(\tau-s)}}\exp\left(-\frac{(x-z)^2}{2\sigma^2(\tau-s)}\right)g(u_n(s,z+\rho (\tau -s))) \mathrm{e}^{-(r + c_M)(\tau-s)}\,\mathrm{d}z\,\mathrm{d}s\right| \\
    &\leq \int_0^\tau\int_{L}^{\infty}\frac{1}{\sqrt{2\pi\sigma^2(\tau-s)}}\exp\left(-\frac{(x-z)^2}{2\sigma^2(\tau-s)}\right)|g(u_n(s,z+\rho (\tau -s)))| \mathrm{e}^{-(r + c_M)(\tau-s)}\,\mathrm{d}z\,\mathrm{d}s \\
    &\leq C^+\int_0^\tau\int_{L}^{\infty}\frac{1}{\sqrt{2\pi\sigma^2(\tau-s)}}\exp\left(-\frac{(x-z)^2}{2\sigma^2(\tau-s)}+az\right)\mathrm{e}^{\overline{\lambda}s-(r + c_M-a\rho)(\tau-s)}\,\mathrm{d}z\,\mathrm{d}s \\
    &+ C^-\int_0^\tau\int_{L}^{\infty}\frac{1}{\sqrt{2\pi\sigma^2(\tau-s)}}\exp\left(-\frac{(x-z)^2}{2\sigma^2(\tau-s)}\right)\mathrm{e}^{\underline{\lambda}s-(r + c_M)(\tau-s)}\,\mathrm{d}z\,\mathrm{d}s \\
    &= C^+\mathrm{e}^{-(r+c_M-a\rho)\tau}\int_0^\tau\mathrm{e}^{(\overline{\lambda}+r+c_M-a\rho)s}\int_{L}^{\infty}\frac{\mathrm{e}^{ax+\frac12\sigma^2(\tau-s)a^2}}{\sqrt{2\pi\sigma^2(\tau-s)}}\exp\left(-\frac{(z-x-\sigma^2(\tau-s)a)^2}{2\sigma^2(\tau-s)}\right)\,\mathrm{d}z\,\mathrm{d}s \\
    &+ C^-\mathrm{e}^{-(r+c_M)\tau}\int_0^\tau\mathrm{e}^{(\underline{\lambda}+r+c_M)s}\int_{L}^{\infty}\frac{1}{\sqrt{2\pi\sigma^2(\tau-s)}}\exp\left(-\frac{(z-x)^2}{2\sigma^2(\tau-s)}\right)\,\mathrm{d}z\,\mathrm{d}s \\
    &= C^+\mathrm{e}^{-(r+c_M-a\rho-\frac12a^2\sigma^2)\tau+ax}\int_0^\tau\mathrm{e}^{(\overline{\lambda}+r+c_M-a\rho-\frac12a^2\sigma^2)s}\frac{1}{\sqrt{\pi}}\int_{\frac{L-x-\sigma^2(\tau-s)a}{\sqrt{2\sigma^2(\tau-s)}}}^{\infty}\exp\left(-y^2\right)\,\mathrm{d}y\,\mathrm{d}s \\
    &+ C^-\mathrm{e}^{-(r+c_M)\tau}\int_0^\tau\mathrm{e}^{(\underline{\lambda}+r+c_M)s}\frac{1}{\sqrt{\pi}}\int_{\frac{L-x}{\sqrt{2\sigma^2(\tau-s)}}}^{\infty}\exp\left(-y^2\right)\,\mathrm{d}y\,\mathrm{d}s \\
    &= \frac12C^+\mathrm{e}^{-(r+c_M-a\rho-\frac12a^2\sigma^2)\tau+ax}\int_0^\tau\mathrm{e}^{(\overline{\lambda}+r+c_M-a\rho-\frac12a^2\sigma^2)s}\operatorname{erfc}\left(\frac{L-x-\sigma^2(\tau-s)a}{\sqrt{2\sigma^2(\tau-s)}}\right)\,\mathrm{d}s \\
    &+ \frac12C^-\mathrm{e}^{-(r+c_M)\tau}\int_0^\tau\mathrm{e}^{(\underline{\lambda}+r+c_M)s}\operatorname{erfc}\left(\frac{L-x}{\sqrt{2\sigma^2(\tau-s)}}\right)\,\mathrm{d}s,
\end{align*}
where we used the well-known complementary error function, \cite[Chap. 2]{Lebedev72},
$$
\operatorname{erfc}\zeta=1-\operatorname{erf}\zeta=\frac{2}{\sqrt{\pi}}\int_{\zeta}^{\infty}\exp(-\xi^2)\,\mathrm{d}\xi.
$$
We analogously obtain
\begin{align*}
    &\left|\int_0^\tau\int_{-\infty}^{-L}\frac{1}{\sqrt{2\pi\sigma^2(\tau-s)}}\exp\left(-\frac{(x-z)^2}{2\sigma^2(\tau-s)}\right)g(u_n(s,z+\rho (\tau -s))) \mathrm{e}^{-(r + c_M)(\tau-s)}\,\mathrm{d}z\,\mathrm{d}s\right| \\
    &\leq C^+\mathrm{e}^{-(r+c_M-a\rho-\frac12a^2\sigma^2)\tau+ax}\int_0^\tau\mathrm{e}^{(\overline{\lambda}+r+c_M-a\rho-\frac12a^2\sigma^2)s}\frac{1}{\sqrt{\pi}}\int_{-\infty}^{\frac{-L-x-\sigma^2(\tau-s)a}{\sqrt{2\sigma^2(\tau-s)}}}\exp\left(-y^2\right)\,\mathrm{d}y\,\mathrm{d}s \\
    &+ C^-\mathrm{e}^{-(r+c_M)\tau}\int_0^\tau\mathrm{e}^{(\underline{\lambda}+r+c_M)s}\frac{1}{\sqrt{\pi}}\int_{-\infty}^{\frac{-L-x}{\sqrt{2\sigma^2(\tau-s)}}}\exp\left(-y^2\right)\,\mathrm{d}y\,\mathrm{d}s \\
    &= \frac12C^+\mathrm{e}^{-(r+c_M-a\rho-\frac12a^2\sigma^2)\tau+ax}\int_0^\tau\mathrm{e}^{(\overline{\lambda}+r+c_M-a\rho-\frac12a^2\sigma^2)s}\operatorname{erfc}\left(\frac{L+x+\sigma^2(\tau-s)a}{\sqrt{2\sigma^2(\tau-s)}}\right)\,\mathrm{d}s \\
    &+\frac12C^-\mathrm{e}^{-(r+c_M)\tau}\int_0^\tau\mathrm{e}^{(\underline{\lambda}+r+c_M)s}\operatorname{erfc}\left(\frac{L+x}{\sqrt{2\sigma^2(\tau-s)}}\right)\,\mathrm{d}s.
\end{align*}
Here, we used that
$$
\frac{2}{\sqrt{\pi}}\int_{-\infty}^{\zeta}\exp(-\xi^2)\,\mathrm{d}\xi
= \operatorname{erf}(\zeta)+1 = 1-\operatorname{erf}(-\zeta) = \operatorname{erfc}(-\zeta).
$$

The first integral in \eqref{e:iter}, that also appears in the solution formula \eqref{e:lin} for the linear case, can be represented in terms of the complementary error function for all three considered contracts. For each of these contracts, the integral has the form
\begin{equation*}
    I = I_1 - I_2 = \int_{\tilde{L}}^{\infty}\frac{1}{\sqrt{2\pi\sigma^2\tau}}\exp\left(-\frac{(x+\rho\tau-y)^2}{2\sigma^2\tau}\right)\mathrm{e}^y\,\mathrm{d}y - \tilde{K}\int_{\tilde{L}}^{\infty}\frac{1}{\sqrt{2\pi\sigma^2\tau}}\exp\left(-\frac{(x+\rho\tau-y)^2}{2\sigma^2\tau}\right)\,\mathrm{d}y
\end{equation*}
with either $\tilde{K}=K$ or $\tilde{K}=K_s$, and $\tilde{L}=-\infty$, $\tilde{L}=\ln K$, or $\tilde{L}=\ln K_t$, respectively. We derive
\begin{align*}
    I_1 &= \int_{\tilde{L}}^{\infty}\frac{1}{\sqrt{2\pi\sigma^2\tau}}\exp\left(-\frac{(y-x-\rho\tau)^2}{2\sigma^2\tau}+y\right)\,\mathrm{d}y \\
    &= \mathrm{e}^{x+(\rho+\frac12\sigma^2)\tau}\int_{\tilde{L}}^{\infty}\frac{1}{\sqrt{2\pi\sigma^2\tau}}\exp\left(-\frac{(y-x-\rho\tau-\sigma^2\tau)^2}{2\sigma^2\tau}\right)\,\mathrm{d}y \\
    &= \frac12\mathrm{e}^{x+(q_s-\gamma_S)\tau}\operatorname{erfc}\left(\frac{\tilde{L}-x-(\rho+\sigma^2)\tau}{\sqrt{2\sigma^2\tau}}\right)
\end{align*}
and
\begin{equation*}
    I_2 = \tilde{K}\int_{\tilde{L}}^{\infty}\frac{1}{\sqrt{2\pi\sigma^2\tau}}\exp\left(-\frac{(y-x-\rho\tau)^2}{2\sigma^2\tau}\right)\,\mathrm{d}y
    = \frac{\tilde{K}}{2}\operatorname{erfc}\left(\frac{\tilde{L}-x-\rho\tau}{\sqrt{2\sigma^2\tau}}\right).
\end{equation*}
For the forward, we have $\tilde{L}=-\infty$ and obtain
\begin{equation*}
    I_{\mathrm{fwd}} = \mathrm{e}^{x+(q_S-\gamma_S)\tau} - K.
\end{equation*}
Furthermore, we have
\begin{equation*}
    I_{\mathrm{call}} = \frac12\mathrm{e}^{x+(q_S-\gamma_S)\tau}\operatorname{erfc}\left(\frac{\ln K-x-(\rho+\sigma^2)\tau}{\sqrt{2\sigma^2\tau}}\right) - \frac12K\operatorname{erfc}\left(\frac{\ln K-x-\rho\tau}{\sqrt{2\sigma^2\tau}}\right)
\end{equation*}
for the call option and
\begin{equation*}
    I_{\mathrm{put}}
    = \frac12\mathrm{e}^{x+(q_S-\gamma_S)\tau}\operatorname{erfc}\left(\frac{\ln K_t-x-(\rho+\sigma^2)\tau}{\sqrt{2\sigma^2\tau}}\right) - \frac12K_s\operatorname{erfc}\left(\frac{\ln K_t-x-\rho\tau}{\sqrt{2\sigma^2\tau}}\right)
\end{equation*}
for the gap option.

\subsection{Numerical implementation of the iteration scheme}\label{ssec:num}

In order to be able to numerically evaluate the integral~\eqref{e:iter} we want to utilize the Gaussian quadrature rule. To achieve this, we need to overcome two major obstacles both resulting from the fact that the values of the $n$-th iteration $u_n$ can be computed only for a finite number of domain points. Obviously, we want the number of computations to be minimal. 

Firstly, the argument of the second improper integral in~\eqref{e:iter} can be evaluated only on a bounded domain. To handle this, we truncate the integration domain while estimating the truncation error as shown in Section \ref{ssec:sc}. The only inconvenience of this approach is that the domain of $u_n$ shrinks with each iteration. Indeed, in order to find the value of $u_{n+1}(x,t)$ one needs to know values of $u_n$ on the rectangle $[0,t] \times \mathcal{I}$ in which $\mathcal{I}$ is a certain neighborhood of the point $x$. This is a consequence of the convolution integral in the iteration step~\eqref{e:iter}.
Thus, we a-priori define the domain in the $(\tau,x)$ space for which we want to compute the values of the derivative contract together with the number of iterations. Then based on the desired truncation error we evaluate backwards the domains for all iterations.

Secondly, due to the shift in the space variable of $u_n$ we either have to cautiously match the discretization points of the consecutive iterations in order to be able to evaluate the integrals or, to circumvent this, we can compute the $n$-th iteration for given rectangular mesh points and subsequently approximate these points via a B-spline basis. The interpolation function is then used as the sub- or supersolution for the next iteration. We choose the order of the basis depending on the option type. When computing forward contract it is more convenient to use a higher-order basis in order to obtain a better interpolation of the function values between the computed points. However, the European call contract and the gap option have non-smooth (even discontinuous) payoff at maturity. Hence, it should be more beneficial to use first-order interpolation to prevent oscillations for such contracts.

\section{Results}\label{sec:results}

To illustrate the scheme, we run numerical computations for the European call option, the forward contract and the gap option in the software Wolfram Mathematica. The Gauss quadrature was coded from scratch without using the internal function \texttt{NIntegrate[]}, the interpolation of the computed points was made by the \texttt{Interpolation[]} function. The model parameters were almost exactly the same for all three types of the contract, see Table~\ref{t:parameters}. 
\begin{table}[ht]
    \centering
    \begin{tabular}{c|l|l}
         Parameter & Value & Note \\
         \hline
         $q_S$ & $0.1$ &\\
         $\gamma_S$ & $0.02$& \\
         $\sigma$ & $0.5$&  \\
         $r$ & $0.1$ &\\
         $\lambda_B$ & $0.02$ &\\
         $\lambda_C$ & $0.03$& \\
         $R_B$ & $0.7$ &\\
         $R_C$ & $0.8$ &\\ 
         $s_F$ & $0.05$& \\
         $K$ & $15$ & for the European call and the forward contract\\
         $K_s$ & $15$ & for the gap option \\
         $K_t$ & $12$ & for the gap option \\
    \end{tabular}
    \caption{Choice of parameters for the model~\eqref{e:check-v} and various contracts.}
    \label{t:parameters}
\end{table}
We considered the time to maturity to be $T=2$, i.e., $\tau \in [0,2]$ with the time discretization $\delta_\tau = 0.02$. The final size of the price domain was set to be $S\in [K-5, K+5] = [10,20]$, i.e., $x \in [x_\mathrm{min},x_\mathrm{max}] := [\log(K-5),\log(K+5)]$ ($x \in [\log(K_t-5),\log(K_s+5)]$, respectively, for the gap option). The truncation error estimates in Section~\ref{ssec:sc} show that the corresponding error (while restricting the improper integral in~\eqref{e:iter} to a finite interval $[-L,L]$) is almost solely dependent on the difference $|L-x|$. Note that
\[
\int_0^2 \mathrm{erfc}\left( \frac{|L-x|}{\sqrt{2-s}} \right)=10^{-16}
\]
has solution $|L-x| \approx 8$ in which the integral underflows the machine precision. Thus, we chose to truncate the price domain in each step by $\epsilon=10$. Computing $5$ iterations and taking into account that for $\rho = -0.42<0$ the $n$-th iteration is computed over the price domain $x_n \in [x_\mathrm{min}-(3-n)(\epsilon+\rho T),x_\mathrm{max}+(3-n)\epsilon]$, $n=1, \ldots, 5$. If $\rho$ were to be positive the term $\rho T$ would appear in the upper bound of the interval, see~\eqref{e:iter}. The discretization step in the price was $\delta_x =0.02$. 

\begin{figure}[ht]
    \begin{subfigure}{0.45\textwidth}
    \includegraphics[width=\textwidth]{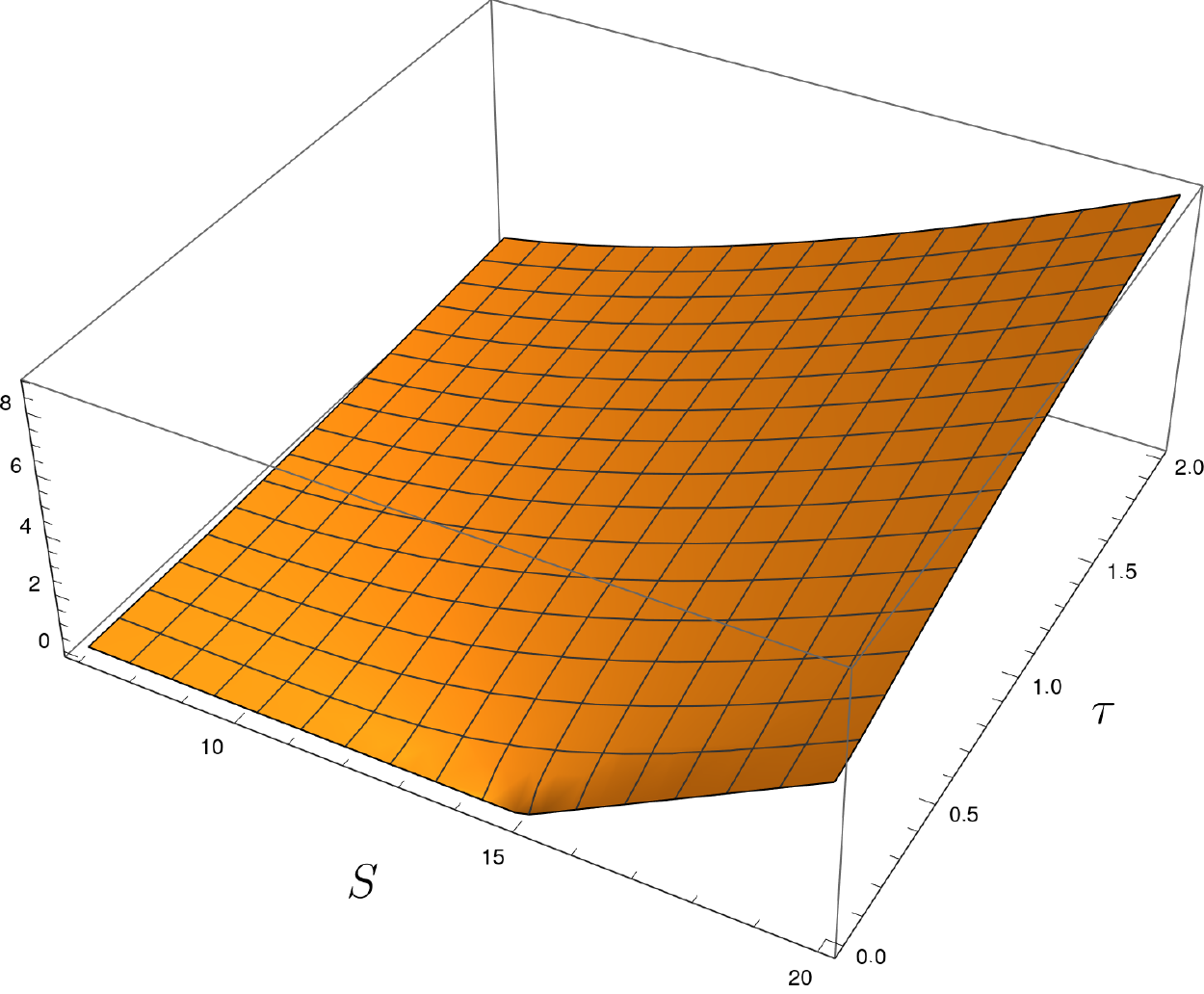}
    \caption{}
    \label{f:call-surf}
\end{subfigure}
\hfill
\begin{subfigure}{0.45\textwidth}
    \includegraphics[width=\textwidth]{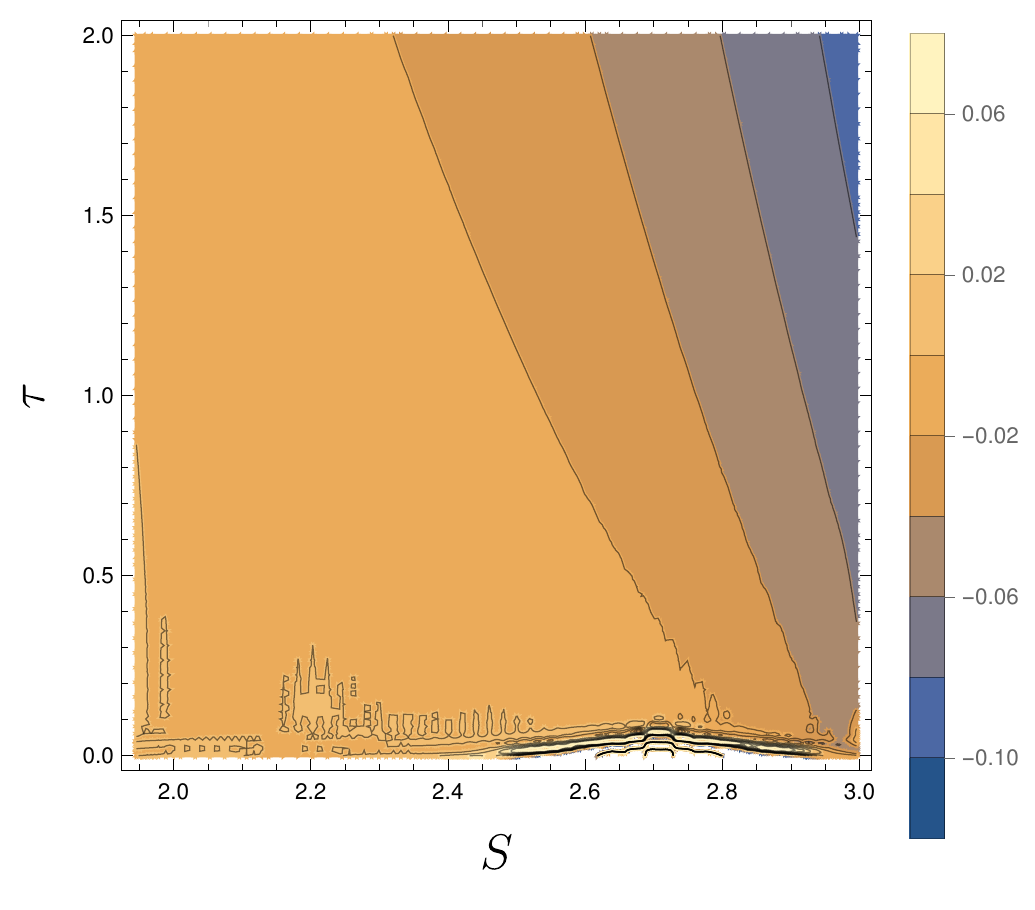}
    \caption{}
    \label{f:call-res}
\end{subfigure}

    \caption{The first panel (a) depicts the resulting surface of the European call option after $5$ iterations. The residual $\mathcal{D}$ is depicted in the right panel (b). }
    \label{f:call}
\end{figure}

For the reference, we also include the results of the numerical computations for the European call option although the price cannot be negative due to the weak maximum principle. See Figure~\ref{f:call}. Thanks to the nature of our approach, we have the piecewise polynomial (B-spline) interpolation of the resulting surface at hand and thus we can compute the residual of our super solution defined by~\eqref{e:residual}. 
The sequence of sub- (super-) solutions preserves the ordering $\underline{u_{n}} \le \overline{u_m}$ for any combination of iterations $m,n \in \mathbb{N}$. Although the sign of the residual should also be preserved, i.e., $\mathcal{D}(\underline{u_{n}}) \le 0$, $\mathcal{D}(\overline{u_{n}}) \ge 0$, this is only true for some couples of initial iterations. When the iterations are close enough to the solution, the truncation and the discretization error disrupt the ordering of the residuals. 

\begin{figure}[ht]
    \begin{subfigure}{0.45\textwidth}
    \includegraphics[width=\textwidth]{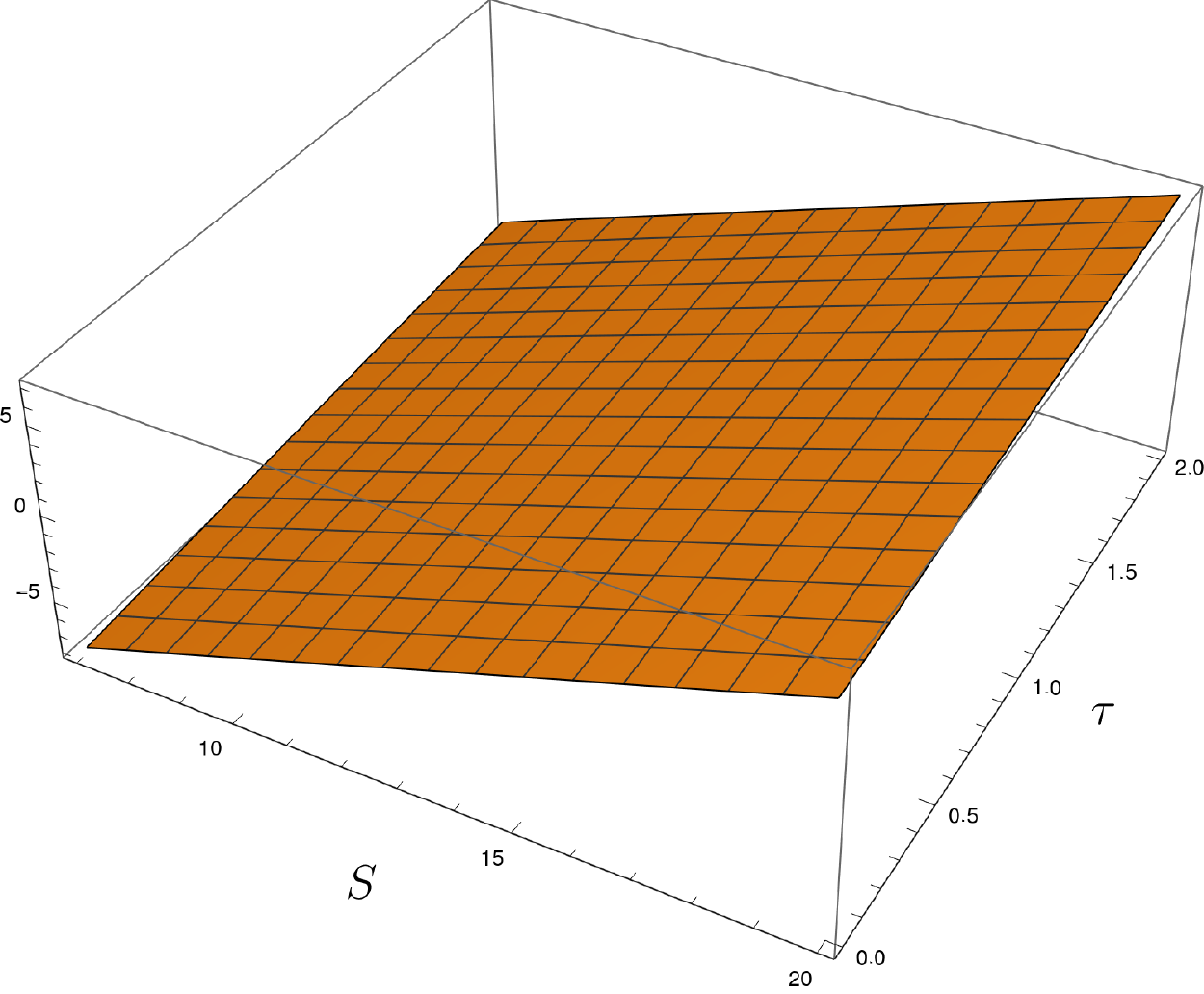}
    \caption{}
    \label{f:fwd-surf}
\end{subfigure}
\hfill
\begin{subfigure}{0.45\textwidth}
    \includegraphics[width=\textwidth]{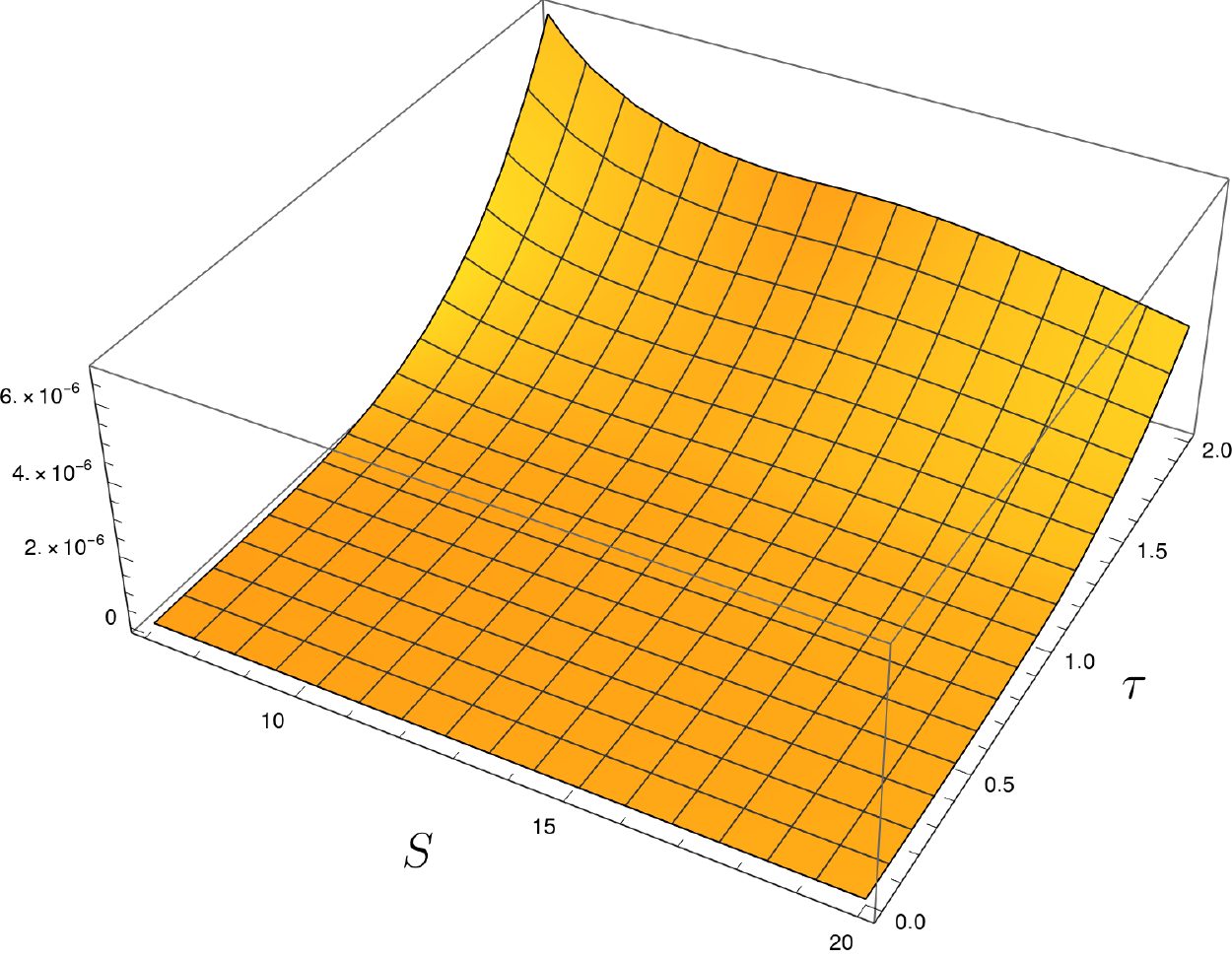}
    \caption{}
    \label{f:fwd-res}
\end{subfigure}
    
    \caption{The first panel (a) depicts the resulting surface of the forward contract after $5$ iterations. The difference between the final sub- and supersolutions are plotted in the right panel (b).}
    \label{f:fwd}
\end{figure}

The results for the forward contract are depicted in Figure~\ref{f:fwd}. 

\begin{figure}[ht]
    \begin{subfigure}{0.45\textwidth}
    \includegraphics[width=\textwidth]{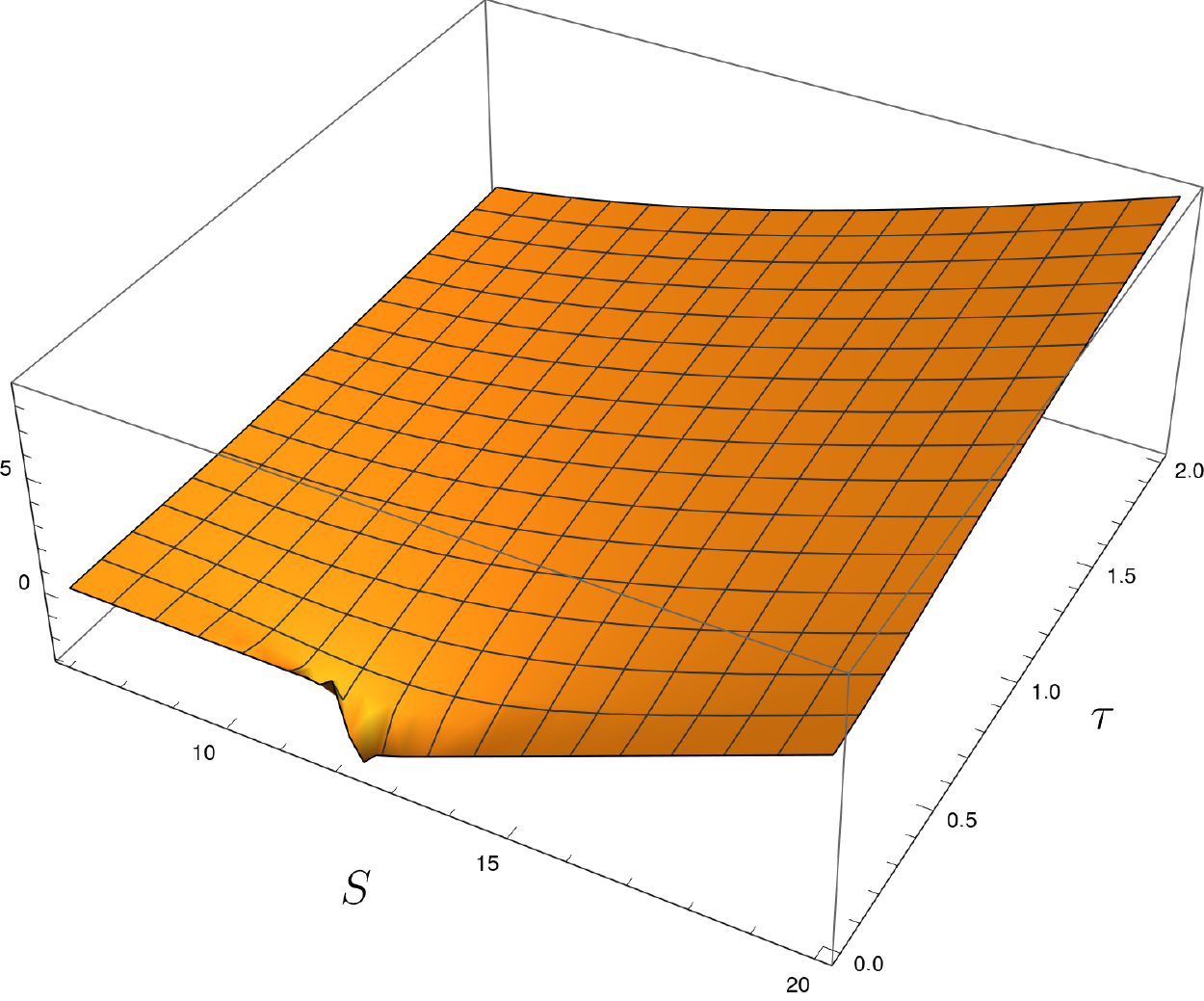}
    \caption{}
    \label{f:gap-surf}
\end{subfigure}
\hfill
\begin{subfigure}{0.45\textwidth}
    \includegraphics[width=\textwidth]{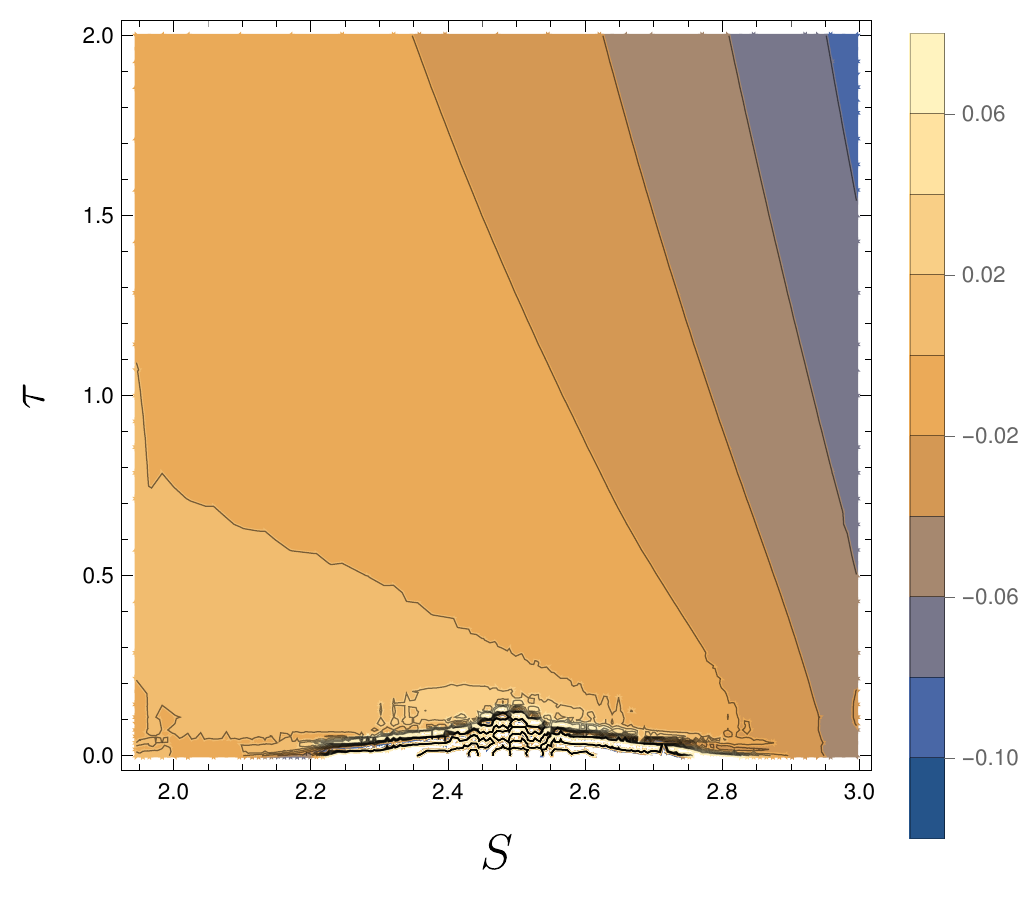}
    \caption{}
    \label{f:gap-res}
\end{subfigure}

    \caption{The first panel (a) depicts the resulting surface of the gap option after $5$ iterations. The residual~\eqref{e:residual} is depicted in the right panel (b).}
    \label{f:gap}
\end{figure}

The behavior of the residual for the gap option in Figure~\ref{f:gap} shows similar numerical properties as of the European call options. It also illustrates the difficulty of approximating a nonsmooth (call option) or even discontinuous function (gap option) with a smooth basis.

\section{Conclusion}\label{sec:conclusion}

In this article, we studied a nonlinear PDE that appears in counterparty risk models provided that the mark-to-market value is chosen to be the adjusted price of the contract. We introduced a monotone iteration scheme using sub- and supersolutions based on recent research in \cite{AlziaryTakac22}. Semi-closed formulas for each iteration step were established together with estimates on the truncation error for the improper integral term. We showed that solutions for contracts with non-negative payoff like the call option are given by a simple modification of the riskless price and, thus, applied the iteration scheme to the forward and the gap call option which are contracts with sign-changing payoff. 

The numerical computations we provided serve as a proof-of-concept. 
Their purpose is to show viability of the monotone iterations approach. 
Indeed, one does not need excessive number of iterations to obtain reasonably precise solution, see, e.g., Figure~\ref{f:call-res}. 
The whole scheme can in addition be viewed as explicit without the need of employing sophisticated equation solvers. 
There is also one specific property of the scheme; if the volatility $\sigma$ is low, then we end up with advection-diffusion with large P\'{e}clet number $2(q_S-\gamma_S-\tfrac{1}{2} \sigma^2)/\sigma^2$ ; naturally assuming $q_S \ne \gamma_S$. 
This is a situation in which numerical schemes can in general lose stability,~\cite{Roos2008robust}. 
In contrary, the integrands in~\eqref{e:iter} have thinner tails which enables us to further truncate the integral domain in our setting, see also~\ref{ssec:sc}. This decreases the computational demands. 

There is also some room for improvement of the numerical implementation.
 To name just few options: the use of programming languages more suitable for numerical computations such as Fortran or C++, the incorporation of adaptive spatial discretization reflexing the shrinking spatial region, the optimization of the order of the Gaussian quadrature, and the optimization of the interpolation order. 
 One also does not need to compute both upper and lower solution and sacrifice the control of the scheme for doubling the speed. 

There are very few downsides of the approach.  
 As stated before, one has to estimate the number of iterations which is sufficient to obtain a solution in advance; before setting the iteration scheme. 
 Moreover, our monotone numerical iteration scheme suffer from the cases of ``overshooting'' in which the solution residuals lose proper ordering. We were not able to find any reasonable reference in literature. Authors usually use the proximity of the lower and the upper solutions as the measure of accuracy of the scheme. It is not common to examine the solution residual $\mathcal{D}$.  
 We argue that slowing the iteration scheme might mitigate the issue by introducing a parameter $\alpha \in (0,1]$ which interpolates between the current and the next iteration 
 \begin{equation}\label{e:slow-iteration}
     u_{n+1}:= \alpha \mathcal{F}(u_n)+ (1-\alpha) u_n
 \end{equation} in which $\mathcal{F}$ is a one-step iteration operator as in~\eqref{e:iter}. Our scheme~\eqref{e:iter} is a special case of~\eqref{e:slow-iteration} with $\alpha =1$. On the other hand, this would inherently increase the number of iterations. 

To solve the sequence of PDEs \eqref{e:sub} and \eqref{e:super}, traditional numerical methods such as finite differences or finite elements methods could be used if monotonicity preserving modifications are applied, see for example the discussion for finite differences method in \cite{Strikwerda2004} and for finite element methods in \cite{Thomee2006}. Comparison with these types of methods was beyond the scope of presented paper and opens several interesting directions in further research. Other improvements might include variable order of the interpolation, together with adaptive discretization, etc.


{\small

}

\end{document}